\documentclass[12pt]{amsart}
\usepackage{amsfonts,amssymb,amsmath}
\usepackage[english]{babel}
\usepackage{amsfonts,amssymb}
\usepackage{geometry}
\newtheorem{theorem}{Theorem}
\newtheorem{lemma}{Lemma}
\newenvironment{definition}[1][Definition]{\begin{trivlist}
\item[\hskip \labelsep {\bfseries #1}]}{\end{trivlist}}

\begin{document}

\author{Andrii Goriunov, Vladimir Mikhailets}

\address{Institute of Mathematics of National Academy of Sciences of Ukraine \\
         Tereshchenkivska str., 3 \\
         Kyiv-4 \\
         Ukraine \\
         01601}

\email[Andrii Goriunov]{goriunov@imath.kiev.ua}         
\email[Vladimir Mikhailets]{mikhailets@imath.kiev.ua}

\title[Resolvent convergence of Sturm-Liouville operators]{Resolvent convergence of Sturm-Liouville operators with singular potentials}
\begin{abstract}
In this paper we consider the Sturm-Liuoville operator in the Hilbert space $L_2$ with the singular complex
potential of $W^{-1}_2$ and two-point boundary conditions.
For this operator we give sufficient conditions for norm resolvent approximation by the operators
of the same class.
\end{abstract}


\maketitle

\section{Main result}

Let on a compact interval $[a,b]$ the formal differential expression
\begin{equation}\label{vyraz}
l(y) = -y''(t) + q'(t)y(t), \qquad q(\cdot) \in L_2([a,b],
\mathbb{C}) =: L_2.
\end{equation} be given.

This expression can be defined as the Shin-Zettl \cite{EM} quasi-differential expression
with following quasi-derivatives \cite{S-Sh}:

$$
 D^{[0]} y = y, \quad D^{[1]} y = y' - qy, \quad  D^{[2]} y =
-(D^{[1]} y)' - qD^{[1]} y - q^2y.
$$

In this paper we consider the set of quasi-differential expressions $l_\varepsilon(\cdot)$
of the form (\ref{vyraz}) with potentials
$q_\varepsilon(\cdot) \in L_2,$ $\varepsilon \in [0, \varepsilon_0]$.
In the Hilbert space $L_2$ with norm $\|\cdot\|_2$
each of these expressions generates a dense closed quasi-differential operator
$L_\varepsilon y :=
l_\varepsilon(y),$

 $$Dom(L_\varepsilon) := \{y
\in L_2: \exists  D_\varepsilon^{[2]} y \in L_2;\quad
\alpha(\varepsilon)\mathcal{Y}_a(\varepsilon)+\beta(\varepsilon)\mathcal{Y}_b(\varepsilon)=0\},$$

\noindent where matrices $\alpha(\varepsilon),\beta(\varepsilon)\in
\mathbb{C}^{2\times 2},$ and vectors $$\mathcal{Y}_a(\varepsilon):=\{
y(a),D^{[1]}_\varepsilon y(a)\},\quad \mathcal{Y}_b(\varepsilon):=\{
y(b),D^{[1]}_\varepsilon y(b)\} \in \mathbb{C}^2.$$

Recall that operators $L_\varepsilon$ converge to $L_0$
in the sense of norm resolvent convergence,
$L_\varepsilon \stackrel{R}{\rightarrow} L_0$, if
there exists a number  $\mu \in \mathbb{C}$
such that $\mu \in \rho(L_0)$ and
$\mu \in \rho(L_\varepsilon)$ (for all sufficiently small $\varepsilon$) and
$$\|(L_\varepsilon - \mu)^{-1} - (L_0 - \mu)^{-1}\| \rightarrow 0, \quad \varepsilon \rightarrow +0.$$

This definition does not depend on the point $\mu \in \rho(L_0)$
\cite{K}.

For the case where matrices
$\alpha(\varepsilon),\beta(\varepsilon)$ do not depend on
$\varepsilon$, paper \cite{S-Sh} gives following

\begin{theorem} \label{S-Sh}

Suppose $\|q_\varepsilon - q_0\|_2 \rightarrow 0$ for $\varepsilon
\rightarrow +0$ and the resolvent set of the operator $L_0$ is not empty.
Then $L_\varepsilon \stackrel{R}{\rightarrow} L_0$.
\end{theorem}

Our goal is to generalize Theorem \ref{S-Sh} onto the case of boundary conditions
depending on $\varepsilon$ and to weaken conditions on potentials applying results of papers \cite{M2,M}.

Denote by $c^{\vee}(t) := \int_a^tc(x)dx$ and by $\|\cdot\|_C$ the sup-norm.

\begin{theorem} \label{G-M}
Suppose the resolvent set of the operator $L_0$ is not empty and
for $\varepsilon \rightarrow +0$:

$1) \quad \|q_\varepsilon\|_2 = O (1)$;

$2) \quad \|(q_\varepsilon - q_0)^\vee\|_C \rightarrow 0$;

$3) \quad \|(q^2_\varepsilon - q^2_0)^\vee\|_C \rightarrow 0$;

$4) \quad \alpha(\varepsilon)\longrightarrow \alpha(0),\quad
\beta(\varepsilon)\longrightarrow \beta(0).$

Then $L_\varepsilon \stackrel{R}{\rightarrow} L_0$.

\end{theorem}

Note that condition 3) is not additive.

Condition 1) (taking into account 2), 3)) may be weakened in several directions.

Actually we will prove a stronger statement on the considered operators'
Green functions' convergence with respect to the norm $\|\cdot\|_\infty$ of the space
$L_\infty$ on the square $[a, b] \times [a, b]$ .

\section{Comparison of Theorems \ref{S-Sh} and \ref{G-M}}

We are going to show that if
$\|q_\varepsilon - q_0\|_2 \rightarrow 0$, $\varepsilon \rightarrow +0,$
then conditions 1), 2), 3) of Theorem \ref{G-M} are true.

Indeed,
$\|q_\varepsilon\|_2 \leq \|q_\varepsilon - q_0\|_2 + \|q_0\|_2 = O (1)$.

Also
\begin{align*}
 |\int_a^t(q_\varepsilon - q_0)ds| &\leq
\int_a^b|q_\varepsilon - q_0|ds \leq (\int_a^b|q_\varepsilon -
q_0|^2 ds)^{1/2}(b - a)^{1/2} \rightarrow 0, \, \varepsilon
\rightarrow +0.\\
 |\int_a^t(q_\varepsilon^2 -
q_0^2)ds| &\leq \int_a^b|q_\varepsilon^2 - q_0^2|ds \leq
{\int_a^b|q_\varepsilon - q_0| |q_\varepsilon + q_0|ds} \leq\\
&\leq (\int_a^b|q_\varepsilon -
q_0|^2ds)^{1/2}(\int_a^b|q_\varepsilon + q_0|^2ds)^{1/2} \rightarrow
0, \quad \varepsilon \rightarrow +0.
\end{align*}
Following example proves Theorem \ref{G-M} to be stronger than Theorem \ref{S-Sh}.
\vskip 0.2 cm
\textsc{Example 1.} Suppose $q_0(t) \equiv 0$, $q_\varepsilon(t) =
e^{it/\varepsilon}$, $t \in [0, 1]$.

The set of operators $L_\varepsilon$ defined by these potentials does not
satisfy assumptions of Theorem \ref{S-Sh} because
$$\|q_\varepsilon - q_0 \|^2_2 = \|q_\varepsilon\|^2_2 =
\int_0^1|q_\varepsilon|^2ds \equiv 1.$$

It is evident that functions $q_\varepsilon(\cdot)$ do not converge to 0
even with respect to the Lebesgue measure. However, they satisfy conditions
1), 2), 3) of Theorem \ref{G-M}.
Indeed, $\|q_\varepsilon\|_2 \leq 1$.
Moreover,
\begin{align*}
&\|q^\vee_\varepsilon\|_C =  \|\int_0^t e^{is/\varepsilon}ds\|_C
\leq
2\varepsilon \rightarrow 0,\quad \varepsilon \rightarrow +0.\\
&\|(q^2_\varepsilon)^\vee\|_C = \|\int_0^t
(e^{is/\varepsilon})^2ds\|_C \leq \varepsilon \rightarrow 0, \quad
\varepsilon \rightarrow +0.
\end{align*}

\section{Preliminary result}

Consider a boundary-value problem

$$
    y'(t;\varepsilon)=A(t;\varepsilon)y(t;\varepsilon)+f(t;\varepsilon),\quad
t \in[a,b], \quad \varepsilon \in [0, \varepsilon_0] \eqno
(3.1_\varepsilon)
$$
$$
    U_{\varepsilon}y(\cdot;\varepsilon)= 0, \eqno
(3.2_\varepsilon)
$$
\addtocounter{equation}{2}
where matrix functions $A(\cdot,\varepsilon) \in  L_{1}^{m\times m}$,
vector-functions $f(\cdot,\varepsilon) \in L_{1}^{m}$,\quad
and linear continuous operators
$U_{\varepsilon}:C([a,b];\mathbb{C}^{m})\rightarrow\mathbb{C}^{m}.$

We recall from \cite{M2, M}
\begin{definition}\label{M^m}
Denote by
$\mathcal{M}^{m}[a,b]=:\mathcal{M}^{m},$ $ m\in \mathbb{N}$
the class of matrix functions
$R(\cdot;\varepsilon):[0,\varepsilon_0]\rightarrow L_1 ^{m\times m}$,
such that the solution of the Cauchy problem
$$ Z'(t;\varepsilon)= R( t;\varepsilon)Z(t;\varepsilon), \quad Z(a;\varepsilon) = I_m$$
satisfies the limit condition
$$\lim\limits_{\varepsilon \rightarrow +0} \|Z(\cdot;\varepsilon) - I_m\|_C =0.$$
\end{definition}
Sufficient conditions for $R(\cdot;\varepsilon)\in \mathcal{M}^{m}$
derive from \cite{Levin}.
To prove Theorem \ref{G-M} we apply the simplest of them
$$\|R(\cdot;\varepsilon)\|_1 = O (1), \quad \|R^{\vee}(\cdot;\varepsilon)\|_C \rightarrow 0,$$
where $\|\cdot\|_1$ is the norm in  $L_1^{m \times m}$.

Paper \cite{M} gives the following general

\begin{theorem}
\label{1 limit G} Suppose \hangindent=0.42cm
\begin{align*}
 1) \quad &\text{the homogeneous limit boundary-value problem~} (3.1_0),(3.2_0) \text{~with~} f(\cdot;0) \equiv 0  \\
 &\text{~has only zero solution;}\\
  2)\quad &A(\cdot;\varepsilon)-A(\cdot;0)\in \mathcal{M}^m;\\
  3)\quad &\|U_{\varepsilon} - U_{0}\|\rightarrow 0,\quad
\varepsilon\rightarrow +0.
\end{align*}

Then for sufficiently small $\varepsilon$
Green matrices $G(t, s; \varepsilon)$ of problems $(3.1_\varepsilon)$,
 $(3.2_\varepsilon)$ exist
and on the square $[a,b]\times [a,b]$
\begin{equation}\label{G}
 \|G(\cdot,\cdot;\varepsilon)-G(\cdot,\cdot;0)\|_\infty \rightarrow 0,\quad
\varepsilon\rightarrow +0.
\end{equation}
\end{theorem}

Condition 3) of Theorem \ref{1 limit G} cannot be replaced by a weaker condition
of the strong convergence of the operators $U_{\varepsilon} \stackrel{s}\rightarrow U_{0}$ \cite{M}.
However, one may easily see that for multi-point "boundary" operators
$$U_\varepsilon y := \sum\limits_{k=1}^n B_k(\varepsilon) y(t_k), \quad \{t_k\} \subset [a, b],
\quad B_k(\varepsilon) \in \mathbb{C}^{m\times m},\quad n \in
\mathbb{N},$$
both conditions of strong and norm convergence are equivalent to
 $${\|B_k(\varepsilon) - B_k(0)\| \rightarrow
0,}\quad \varepsilon\rightarrow +0, \quad k \in \{1, ..., n\}.$$

\section{Proof of Theorem \ref{G-M}}

We give two lemmas to apply Theorem \ref{1 limit G} to proof of Theorem \ref{G-M}.

\begin{lemma}\label{lemm1}
Function $y(t)$ is a solution of a boundary-value problem
\begin{equation}\label{D^2}
 D^{[2]}_\varepsilon y(t)= f(t;\varepsilon) \in L_2 ,\quad\varepsilon\in [0,\varepsilon_0],
\end{equation}
\begin{equation}\label{alpha+beta}
  \alpha(\varepsilon)\mathcal{Y}_a(\varepsilon)+
  \beta(\varepsilon)\mathcal{Y}_b(\varepsilon)=0.
\end{equation}
if and only if vector-function $w(t) = (y(t),
D^{[1]}_\varepsilon y(t))$ is a solution of a boundary-value problem
\begin{equation}\label{diff eq}
w'(t)=A(t;\varepsilon)w(t) + \varphi(t;\varepsilon),
\end{equation}
\begin{equation}\label{diff alpha+beta}
\alpha(\varepsilon)w(a)+
  \beta(\varepsilon)w(b)=0,
  \end{equation}
where matrix function
\begin{equation}\label{A matrix}
A(\cdot;\varepsilon):=\left ( \begin{array}{cc}
q_\varepsilon& 1\\
-q_\varepsilon^2& -q_\varepsilon
\end{array}\right) \in L_1^{2\times 2},
\end{equation}
and $\varphi(\cdot;\varepsilon) := (0, -f(\cdot;\varepsilon))$.
\end{lemma}

\textsc{Proof.} Consider the system of equations
 $$\left\{
\begin{array}{l}
  ( D^{[0]}_\varepsilon y(t))' = q_\varepsilon(t)D^{[0]}_\varepsilon y(t) + D^{[1]}_\varepsilon y(t)    \\
   ( D^{[1]}_\varepsilon y(t))' = - q_\varepsilon^2(t)D^{[0]}_\varepsilon y(t) - q_\varepsilon(t)D^{[1]}_\varepsilon y(t) - f(t; \varepsilon) \\
\end{array}
\right.$$

If $y(\cdot)$ is a solution of equation (\ref{D^2}),
then definition of quasi-derivatives derives that $y(\cdot)$
is a solution of this system.
On the other hand with
$$w(t) = (D^{[0]}_\varepsilon y(t), D^{[1]}_\varepsilon y(t))
\qquad \text{and} \qquad
\varphi(t;\varepsilon) = (0, -f(t;\varepsilon))$$
this system may be rewritten in the form of equation (\ref{diff eq}).

As $\mathcal{Y}_a(\varepsilon) = w(a)$, $\mathcal{Y}_b(\varepsilon) = w(b)$
then it is evident that boundary conditions (\ref{alpha+beta})
are equivalent to boundary conditions (\ref{diff
alpha+beta}).

\begin{lemma}\label{Gamma exist}
Let the assumption
\begin{itemize}
\item [$(\mathcal{E})$]Homogeneous boundary-value problem $D^{[2]}_0 y(t)=0, \quad
\alpha(0)\mathcal{Y}_a(0)+ \beta(0)\mathcal{Y}_b(0) = 0$
has only zero solution
\end{itemize}
be fulfilled.
Then for sufficiently small $\varepsilon$ Green function
$\Gamma(t,s;\varepsilon)$ of the semi-homogeneous boundary problem
(\ref{D^2}), (\ref{alpha+beta}) exists
and $$ \Gamma(t,s;\varepsilon) = -g_{12}(t,s;\varepsilon) \qquad\mbox{a. e.,}$$
\noindent  where $g_{12}(t,s;\varepsilon)$ is the corresponding element of the Green's matrix
$$G(t,s;\varepsilon)=(g_{ij}(t,s;\varepsilon))_{i,j=1}^2 $$
of two-point vector boundary-value problem
(\ref{diff eq}), (\ref{diff alpha+beta}).
\end{lemma}

\textsc{Proof.} Taking into account Theorem \ref{1 limit G} and Lemma \ref{lemm1}
assumption ($\mathcal{E}$) derives that homogeneous boundary-value problem
$$w'(t)=A(t;\varepsilon)w(t), \quad \alpha(\varepsilon)w(a) +
\beta(\varepsilon)w(b)=0$$
for sufficiently small $\varepsilon$ has only zero solution.

Then for problem (\ref{diff eq}), (\ref{diff alpha+beta}) Green matrix
$$G(t,s, \varepsilon)=(g_{ij}(t,s))_{i,j=1}^2\in L_\infty^{2\times 2}$$
exists and the unique solution of (\ref{diff eq}), (\ref{diff alpha+beta})
is written in the form
$$w_\varepsilon(t)=\int\limits_a ^b
G(t,s;\varepsilon)\varphi(s;\varepsilon) ds, \quad t\in [a,b], \quad
\varphi(\cdot;\varepsilon)\in L_2 .$$

The last equality can be written in the form
$$\left\{
\begin{array}{l}
    D^{[0]}_\varepsilon y_\varepsilon(t) = \int\limits_a^b g_{12}(t,s;\varepsilon)(-\varphi(s;\varepsilon))ds \\
    D^{[1]}_\varepsilon y_\varepsilon(t) = \int\limits_a^b g_{22}(t,s;\varepsilon)(-\varphi(s;\varepsilon))ds, \\
\end{array}
\right.$$
where $y_\varepsilon(\cdot)$ is the unique solution of problem
(\ref{D^2}), (\ref{alpha+beta}).
This implies the assertion of Lemma
\ref{Gamma exist}.

Now, passing to the proof of Theorem  \ref{G-M}, we note that since
$$(q_\varepsilon + \mu)^2 - (q_0 + \mu)^2 = (q_\varepsilon^2 - q_0^2)
+ 2\mu(q_\varepsilon - q_0),$$
in view of conditions 2), 3) we can assume without loss of generality that $0 \in \rho(L_0)$.
Let's prove that
$$\sup\limits_{\|f\|_2 =
1} \|L_\varepsilon^{-1}f - L_0^{-1}f\| \rightarrow 0, \quad \varepsilon
\rightarrow +0.$$

Equation $L_\varepsilon^{-1}f = y_\varepsilon$ is equivalent to the relation
$L_\varepsilon y_\varepsilon = f$, that is $y_\varepsilon$ is the solution
of the problem (\ref{D^2}), (\ref{alpha+beta})
and due to inclusion  $0 \in \rho(L_0)$ the assumption ($\mathcal{E}$) of Lemma \ref{Gamma exist} holds.
Conditions 1)--3) of Theorem \ref{G-M}
imply that $A(\cdot;\varepsilon)-A(\cdot;0)\in \mathcal{M}^2$,
where $A(\cdot;\varepsilon)$ is given by (\ref{A matrix}).
Therefore assumption of Theorem \ref{G-M} derives that
assumption of Theorem \ref{1 limit G} for problem (\ref{diff eq}), (\ref{diff alpha+beta})
is fulfilled.
This means that Green matrices $G(t,s;\varepsilon)$
of the problems (\ref{diff eq}), (\ref{diff alpha+beta}) exist and
limit relation (\ref{G}) holds.
Taking into account Lemma \ref{Gamma exist}, this implies the limit equality
$$ \|\Gamma(\cdot,\cdot;\varepsilon)-\Gamma(\cdot,\cdot;0)\|_\infty \rightarrow
0,\quad \varepsilon\rightarrow +0.$$

Then

\parshape 3 0cm 10 cm 2.5 cm 10 cm 5 cm 10 cm \noindent$\|L_\varepsilon^{-1}- L_0^{-1}\| =$ $\sup\limits_{\|f\|_2=1} \|\int_a^b[\Gamma(t,s;\varepsilon) -
\Gamma(t,s;0)]f(s)ds\|_2 \leq $\\
${(b - a)^{1/2}}{\sup\limits_{\|f\|_2=1}\|\int_a^b
|\Gamma(t,s;\varepsilon) -
\Gamma(t,s;0)| |f(s)| ds\|_C} \leq$\\
$(b - a)\|\Gamma(\cdot,\cdot;\varepsilon) - \Gamma(\cdot,\cdot;0)\|_\infty
\rightarrow 0, \quad \varepsilon \rightarrow +0$,

\noindent which implies the assertion of Theorem \ref{G-M}.

\section{Three extensions of Theorem \ref{G-M}}

As was already noted, the assumptions of Theorem 2 may be weakened. Let
$$R(\cdot; \varepsilon) := A(\cdot; \varepsilon) - A(\cdot; 0)$$
where $A(\cdot;\varepsilon)$ is given by (\ref{A matrix}).

\begin{theorem} \label{G-M general}
In the statement of Theorem \ref{G-M}, condition 1) can be replaced by any one of the
following three more general (in view of 2) and 3)) asymptotic conditions as
$\varepsilon \rightarrow +0$:

\emph{(I)~} \,\,$\|R(\cdot;
\varepsilon)R^\vee(\cdot;\varepsilon)\|_1 \rightarrow 0;$

\emph{(II)~} $\|R^\vee(\cdot;\varepsilon)R(\cdot; \varepsilon)\|_1
\rightarrow 0;$

\emph{(III)} $\|R(\cdot; \varepsilon)R^\vee(\cdot;\varepsilon) -
R^\vee(\cdot;\varepsilon)R(\cdot; \varepsilon)\|_1 \rightarrow 0.$
\end{theorem}

\textsc{Proof.} The proof of Theorem \ref{G-M general} is similar to the proof of
Theorem \ref{G-M} with following remark to be made.
Condition 2) of Theorem \ref{1 limit G} holds if (see \cite{Levin})
$\|R^\vee(\cdot;\varepsilon)\|_C \rightarrow 0$
and either the condition $R(\cdot;\varepsilon)\|_1 = O (1)$
(as in Theorem \ref{G-M}), or any of three conditions (I), (II), (III)
of Theorem \ref{G-M general} holds.

Following example shows each part of Theorem \ref{G-M general} to be stronger than Theorem \ref{G-M}.

\textsc{Example 2.} Let $q_0(t) \equiv 0$, $q_\varepsilon(t) =
\rho(\varepsilon)e^{it/\varepsilon},$ $t \in [0, 1].$

One may easily calculate that conditions
$$\rho(\varepsilon)
\uparrow \infty, \quad \varepsilon\rho^3(\varepsilon) \rightarrow 0,
\quad \varepsilon \rightarrow +0, $$
imply assumptions 2), 3) of Theorem \ref{G-M} and any one of assumptions
(I), (II), (III) of Theorem \ref{G-M general}.
But assumption 1) of Theorem \ref{G-M}, does not hold because $\|q_\varepsilon -
q_0\|_2\uparrow \infty$.

For Schr\"{o}dinger operators of the form (\ref{vyraz}) on $\mathbb{R}$ with real-valued periodic potential $q'$,
where $q \in L_2^{loc}$,
self-adjointness and sufficient conditions for norm resolvent convergence were established
in \cite{MM}. For other problems related to those studied in \cite{S-Sh}, see also \cite{GM},
\cite{GM2}.

\end{document}